\documentclass[10pt]{article}
\usepackage{amsmath,amssymb,amsthm,amscd}
\numberwithin{equation}{section}

\newtheorem{prop}{Proposition}[section]
\newtheorem{theorem}[prop]{Theorem}

\newtheorem{remark}[prop]{Remark}

\newtheorem{acknowledgment}[prop]{Acknowledgment}

\def\<{\langle}
\def\>{\rangle}
\def\({\left(}
\def\){\right)}

\def\p{\partial}

\def\Ric{{\rm Ric}}
\def\osc{{\rm osc}}
\def\sup{{\rm sup}}
\def\inf{{\rm inf}}
\def\exp{{\rm exp}}

\begin{document}

\title{Convergence Results for Two K\"ahler-Ricci Flows}

\author{the University of Sydney \\
Zhou Zhang \footnote{Partially supported by Australian Research Council Discovery Project: DP110102654.} \\
} 
\date{February 27, 2014}
\maketitle

\begin{abstract}

In this note, we provide some general discussion on the two main versions in the study of K\"ahler-Ricci flows over closed manifolds, aiming at smooth convergence to the corresponding K\"ahler-Einstein metrics with assumptions on the volume form and Ricci curvature form along the flow. 

\end{abstract}

\section{Introduction}

K\"ahler-Ricci flow is the natural version of Ricci flow, introduced by Hamilton in \cite{ham82}, in the context of complex manifold. The study of it was led to the case of Fano manifold, i.e. closed complex manfold with positive first Chern class, after the work \cite{cao} by H. D. Cao. Major works have been carried out in this direction, for example, \cite{chen-tian, tian-zhu}, and some most recent ones including \cite{tian-zzl, jiang} appearing after the breakthrough regarding the sufficient condition on the existence of K\"ahler-Einstein metric by G. Tian and and also Chen-Sun-Donaldson. On the other hand, in recent years, K\"ahler-Ricci flow in the most general setting has also attracted extensive attentions, especially with the program proposed by Tian in \cite{tian02} as the geometric analysis version of the renowned Minimal Model Program in algebraic geometry. A lot of progress has been made on this topic, for example, in \cite{t-znote, song-tian}. 

In this work, we consider the strong convergence for these two main versions of K\"ahler-Ricci flows. In the Fano case, we orginally made the assumptions on the volume form and the Ricci curvature form, while the one on Ricci form can be removed by properly applying the result in \cite{ses-tian}. For the flow in the most general setting, the assumptions are on the global volume and the Ricci curvature form. Let's state the main results below. 

\begin{theorem}
\label{thm:fano-1}
Consider the K\"ahler-Ricci flow over a Fano manifold, (\ref{eq:krf}). If uniformly for all time, the volume form is bounded away from $0$, then the K\"ahler-Einstein metric exists over this manifold and the flow metric converges to it at time infinity.
\end{theorem}

\begin{theorem}
\label{thm:general-conv}
Consider the K\"ahler-Ricci flow, (\ref{eq:gkrf}). If uniformly in the maximal time interval of existence, Ricci curvature form is bounded from below by a fixed real smooth $(1, 1)$-form and the global volume is bounded away from $0$, then the flow exists forever. Furthermore, the manifold has negative first Chern class and the flow metric converges to the K\"ahler-Einstein metric at time infinity. 
\end{theorem}

\noindent{\bf Notation:} in the following, $C$ stands for a positive constant which might be different at places. We also use other letters to indicate constants for better illustration of the idea. 

\begin{acknowledgment}

The author would like to thank Professor G. Tian for introducing him to this interesting field of research and constant encouragement. We also appreciate the interest by people including Kwok-Kun Kwong, Yalong Shi and Ben Weinkove. The discussion with Yalong Shi stimulated the improvement on the result for the Fano case. 

\end{acknowledgment}

\section{Fano K\"ahler-Ricci Flow}

The traditional K\"ahler-Ricci flow over a Fano manifold is one case of the more general K\"ahler-Ricci flow. Of course, it is fairly justified to focus on this specific flow because of its natural relation to the search of the necessary and sufficient condition of the existence of K\"ahler-Einstein metrics over Fano manifolds. 

Over a Fano manifold $X$ of complex dimension $n\geqslant 2$, i.e. the first Chern class  $c_1(X)$ being K\"ahler, we consider the following geometric evolution:
\begin{equation}
\label{eq:krf}
\frac{\p \widetilde\omega_t}{\p t}=-\Ric(\widetilde\omega_t)+\widetilde\omega_t, ~~~~\widetilde\omega_0=\omega
\end{equation}
with the initial metric $\omega$ representing $c_1(X)$, i.e. $[\omega]=[\Ric(\omega)]=c_1(X)\in H^2(X;\mathbb{R})\cap H^{1, 1}(X;\mathbb{C})$.

This flow is known to exist forever and $[\widetilde\omega_t]=c_1(X)$ for all time, and so we set $\widetilde\omega_t=\omega+\sqrt{-1}\p\bar\p u$. There is a unique smooth volume for $\Omega$ such that $\omega=\Ric(\Omega)$ and $\int_X\Omega=[\omega]^n$. Now (\ref{eq:krf}) can be transformed to the following equivalent scalar evolution equation for the metric potential $u$: 
\begin{equation}
\label{eq:skrf}
\frac{\p u}{\p t}=\log\frac{(\omega+\sqrt{-1}\p\bar\p u)^n}{\Omega}+u, ~~~~u|_{t=0}=0.
\end{equation}

Take $t$-derivative to get
$$\frac{\p}{\p t}\(\frac{\p u}{\p t}\)=\Delta\(\frac{\p u}{\p t}\)+\frac{\p u}{\p t},$$
which can be transformed to 
$$\frac{\p}{\p t}\(\frac{\p u}{\p t}-u\)=\Delta\(\frac{\p u}{\p t}-u\)+n-\<\widetilde\omega_t, \omega\>,$$
where $\<\widetilde\omega_t, \omega\>$ means taking trace of $\omega$ with respect to $\widetilde\omega_t$.

Meanwhile, (\ref{eq:skrf}) can be reformulated in the following way which is more of a complex Monge-Ampere look: 
\begin{equation}
\label{eq:cma}
(\omega+\sqrt{-1}\p\bar\p u)^n=e^{\frac{\p u}{\p t}-u}\Omega.
\end{equation}

In \cite{ses-tian}, together with other higher order controls, it is shown that   
$$\osc_{X\times \{t\}} \frac{\p u}{\p t}\leqslant C$$ 
uniformly for all time, which is an essential building block in the recent study of this flow.  

\subsection{Volume Form Bounded from Above}

If we assume for some constant $p>1$, 
$$\int_X e^{p\(\frac{\p u}{\p t}-u\)}\Omega\leqslant C$$ 
uniformly for all time, applying the original $L^\infty$-estimate by Kolodziej in \cite{koj98} (or \cite{kojnotes}) to (\ref{eq:cma}), we know that uniformly for all time, 
$$\osc_{X\times \{t\}}u\leqslant C.$$

\begin{remark}
The assumption above can certainly be weakened using Kolodziej's original result, and we leave it like this for the simplification of terminology. Of course, it would be the case if $\frac{\p u}{\p t}-u\leqslant C$, i.e. the volume form has a uniform upper bound. 
\end{remark}

Then as $\int_X e^{\frac{\p u}{\p t}-u}\Omega=[\widetilde\omega_t]^n=c_1(X)^n$, we know 
$$\vline\frac{\p u}{\p t}-u\vline\leqslant C.$$

Now recall the following inequality for Parabolic Schwarz Lemma: 
$$\(\frac{\p}{\p t}-\Delta\)\log\<\widetilde\omega_t, \omega\>\leqslant C\<\widetilde\omega_t, \omega\>+C.$$ 
Combining it with 
$$\(\frac{\p}{\p t}-\Delta\)\(\frac{\p u}{\p t}-u\)=n-\<\widetilde\omega_t, \omega\>,$$ 
we arrive at 
$$\(\frac{\p}{\p t}-\Delta\)\(\log\<\widetilde\omega_t, \omega\>+E\(\frac{\p u}{\p t}-u\)\)\leqslant -\<\widetilde\omega_t, \omega\>+C$$
for some large enough positive constant $E$. Applying  Maximum Principle in the standard way, together with the uniform bound of $\frac{\p u}{\p t}-u$, we conclude  
$$\<\widetilde\omega_t, \omega\>\leqslant C.$$
Hence we could have the following uniform bound of $\widetilde\omega_t$ as metric for all time: 
$$\frac{1}{C}\omega\leqslant \widetilde\omega_t\leqslant C\omega.$$

It is now classic to obtain uniform $C^k$ bounds for the normalized $u$, and one can then apply the K-energy argument to achieve a sequential convergence with the limit being the K\"ahler-Einstein metric. By the result in \cite{tian-zhu}, we actually have the flow convergence as $t\to\infty$. We summarize all this in the following proposition.

\begin{prop}
\label{prop:v-u-b}
If for some $p>1$, $\int_X e^{p\(\frac{\p u}{\p t}-u\)}\Omega\leqslant C$ uniformly for all time, then the K\"ahler-Eintein metric exists over $X$ and the flow metric of (\ref{eq:krf}) converges to it as $t\to\infty$.  
\end{prop}

This result is by no means new and has appeared essentially in the works of Tian and Zhu, for example. We merely use it to collect some known techniques in this business. 

\subsection{Volume Form Bounded from Below}

Now we move on to the case when the volume form has a uniform positive lower bound, i.e. for all time, 
$$\frac{\p u}{\p t}-u\geqslant -C.$$
In this case, we also assume that for all time, 
$$\Ric(\widetilde\omega_t)\geqslant \alpha$$ 
for a fixed real smooth $(1, 1)$-form $\alpha$ over $X$. By (\ref{eq:krf}), we then have 
$$\alpha\leqslant \Ric(\widetilde\omega_t)=-\frac{\p \widetilde\omega_t}{\p t}+\widetilde\omega_t=\omega+\sqrt{-1}\p\bar\p\(-\frac{\p u}{\p t}+u\),$$
and so for some large enough $F$,
$$F\omega+\sqrt{-1}\p\bar\p\(-\frac{\p u}{\p t}+u\)>0.$$

Now we need the solution of Berndtsson in the very recent work \cite{bo} on the Openness Conjecture proposed by Demailly and Kollar. The result is as follows.
\begin{theorem}
\label{thm:openness}
Let $u$ be a plurisubharmonic function in the unit ball $B\subset\mathbb{C}^n$ and $u\leqslant D$. Using the standard Euclidean measure $d\lambda$, assume that 
$$\int_Be^{-u}d\lambda<A<\infty.$$
Then we have for the half unit ball $B/2$ and $\epsilon=\frac{\delta(n)}{e^DA}$ with some universial dimensional constant $\delta(n)$,
$$\int_{B/2}e^{-(1+\epsilon)u}d\lambda<C(n)\cdot\(e^{-\epsilon D}A+e^{-(1+\epsilon)D}\),$$
where $C(n)$ is a positive constant only depending on $n$.
\end{theorem}

The above statement is more quantitative comparing with the orginal statement in \cite{bo} and certainly a direct consequence of the argument there. In application, we look at the function $-\frac{\p u}{\p t}+u\in PSH_{F\omega}(X)$ and $-\frac{\p u}{\p t}+u\leqslant C$. The integral bound is clear from global cohomology information: 
$$\int_X e^{-(-\frac{\p u}{\p t}+u)}\Omega=c_1(X)^n<\infty.$$
So the above result by Berndtsson tells us for some $\epsilon>0$
$$\int_X e^{(1+\epsilon)(\frac{\p u}{\p t}-u)}\Omega=\int_X e^{-(1+\epsilon)(-\frac{\p u}{\p t}+u)}\Omega\leqslant C$$
uniformly for all time. Again, by Kolodziej's $L^\infty$-estimate, we know from (\ref{eq:cma}):  
$$\osc_{X\times\{t\}}u\leqslant C$$
uniformly for all time. Now the argument in Subsection 2.1 can be recycled to get the existence of K\"ahler-Einstein metric and the convergence of flow metric to it. Let's summarize all this in the following theorem. 

\begin{theorem}
\label{th:v-ric-conv}
If for all time, $\frac{\p u}{\p t}-u\geqslant -C$ and $\Ric(\widetilde\omega_t)\geqslant\alpha$ for a fixed real smooth $(1, 1)$-form $\alpha$ over $X$, then the K\"ahler-Einstein metric exists over $X$ and the flow metric of (\ref{eq:krf}) converges to it as $t\to\infty$.

\end{theorem}

The assumption on Ricci form can be removed by observing that the result in \cite{ses-tian} quoted before allows us to replace $\frac{\p u}{\p t}$ by a constant (the average for example) up to a control term, which can be absorbed by the volume form. Then the plurisubharmonicity of the metric potential $u$ is enough for the application of Theorem \ref{thm:openness}, together with the positive lower volume form bound. So we have proved Theorem \ref{thm:fano-1}. 

\begin{remark}

As pointed out by Yalong Shi, such implication of the flow convergence and existence of K\"ahler-Einstein metric by the one-sided volume form bound can be viewed as a reminisence of the classic result of Tian and Siu on the more direct continuity method setting to solve the K\"ahler-Einstein equation on a Fano manifold. 

\end{remark}

\section{General K\"ahler-Ricci Flow}

In this section, we turn to the study of the following K\"ahler-Ricci flow over a closed K\"ahler manifold $X$ of complex dimension $n\geqslant 2$:
\begin{equation}
\label{eq:gkrf}
\frac{\p\widetilde\omega_t}{\p t}=-\Ric(\widetilde\omega_t)-\widetilde\omega_t, ~~~~\widetilde\omega_0=\omega_0,
\end{equation}
where the initial metric $\omega_0$ is any K\"ahker metric over $X$. The K\"ahler class of the flow metric is changing in general along the flow, and we prefer this version because the cohomology is automatically controlled under it. To study this evolution equation, we set 
$$\omega_t=-\Ric(\omega_0)+e^{-t}\(\omega_0+\Ric(\omega_0)\)$$
and $\widetilde\omega_t=\omega_t+\sqrt{-1}\p\bar\p u$. The following scalar evolution equation of the metric potential is equivalent to (\ref{eq:gkrf})
\begin{equation}
\label{eq:sgkrf}
\frac{\p u}{\p t}=\log\frac{(\omega_t+\sqrt{-1}\p\bar\p u)^n}{\omega^n_0}-u, ~~~~u|_{t=0}=0.
\end{equation}

In \cite{t-znote}, the {\it Optimal Existence Result} is obtained, i.e. the flow exists as long as the class $[\omega_t]$ remains to be K\"ahler. We denote $[0, T)$ as the maximal time interval of existence, where $T\in (0, \infty]$, depending on the cohomology information. The main result of this section is the following.
\begin{theorem}
\label{th:gkrf-conv}
If uniformly for all $t\in [0, T)$, $\Ric(\widetilde\omega_t)\geqslant\alpha$ for a fixed real smooth $(1, 1)$-form $\alpha$ and $[\omega_t]^n>\epsilon>0$ for some $\epsilon$, then $T=\infty$ and actually $c_1(X)<0$ (implying $X$ being projective) and the flow converges to the K\"ahler-Einstein metric as $t\to\infty$. 
\end{theorem}

This is just Theorem \ref{thm:general-conv}. Let's clarify that the above $\alpha$ and $\epsilon$ are independent of the singularity time $T$. 

\begin{proof}

In light of $\Ric(\widetilde\omega_t)\geqslant\alpha$ and applying a simple ODE calculation, we see $\widetilde\omega_t\leqslant C\omega_0$ uniformly for $t\in [0, T)$. It is also clear that  
$$\alpha\leqslant \Ric(\widetilde\omega_t)=-\frac{\p \widetilde\omega_t}{\p t}-\widetilde\omega_t=\Ric(\omega_0)+\sqrt{-1}\p\bar\p\(-\frac{\p u}{\p t}-u\),$$
and so for some large enough $F$,  
$$F\omega_0+\sqrt{-1}\p\bar\p\(-\frac{\p u}{\p t}-u\)>0.$$

Applying a result in \cite{tian-87} which is the manifold version of a classic result by H\"ormander in \cite{hor}, we have for any $t\in [0, T)$, a constant $a>0$ and $C$ such that 
$$\int_X \exp\(a\(\sup_{X\times \{t\}}\(-\frac{\p u}{\p t}-u\)+\frac{\p u}{\p t}+u\)\)\omega^n_0\leqslant C,$$
where obviously we could make sure  $a\in (0, 1]$. Now we have 
\begin{equation}
\begin{split}
\exp\(a\cdot\inf_{X\times \{t\}}\(\frac{\p u}{\p t}+u\)\)
&\geqslant C\int_X \exp\(a\(\frac{\p u}{\p t}+u\)\)\omega^n_0 \\
&\geqslant C\int_X \exp\(a\(\frac{\p u}{\p t}+u-C\)\)\omega^n_0 \\
&\geqslant C\int_X \exp\(\frac{\p u}{\p t}+u-C\)\omega^n_0 \\
&= C[\widetilde\omega_t]^n=C[\omega_t]^n>\epsilon>0,
\end{split} \nonumber
\end{equation}
where for the third $\geqslant$, we make use of $\frac{\p u}{\p t}+u\leqslant C$ which are known by standard Maximum Principle argument (as in \cite{r-blow-up} for example). Hence we conclude
$$\frac{\p u}{\p t}+u\geqslant -C,$$ 
which gives the uniform volume form lower bound. Together with $\widetilde\omega_t\leqslant C\omega_0$ pointed out at the beginning, we conclude that uniformly for all $t\in [0, T)$,  
$$\frac{1}{C}\omega_0\leqslant \widetilde\omega_t\leqslant C\omega_0,$$
where $C$ does not depend on $T$. 

As in \cite{r-blow-up}, we can then see that there is no finite time singularity and $[\omega_t]$ is K\"ahler also for $t=\infty$, i.e. $-c_1(X)$ is K\"ahler. In this case, the convergence of flow metric to the K\"ahler-Einstein metric is known even when the initial K\"ahler class is not $-c_1(X)$, for example, in \cite{thesis}. The proof is finished. 

\end{proof}

\begin{remark}

The above argument is similar to and also extends that in \cite{ricci-lower}, since for finite time singularity case $\Ric(\widetilde\omega_t)\geqslant -C\widetilde\omega_t$ implies $\Ric((\widetilde\omega_t)\geqslant \alpha$ for some fixed $\alpha$, but not the other way, at least a priori.

\end{remark}

\section{Final Remarks}

In spite of the essential difference in the existing methods and techniques for the study of these two versions of K\"ahler-Ricci flows, there is similarity illustrated by these results and arguments. It would be interesting to create further combination of the ideas, especially through making use of the recent progress in pluripotential theory.

{\it 
\noindent Zhou Zhang \\ 
\noindent Address: Carslaw Building F07, the School of Mathematics and Statistics \\
the University of Sydney, NSW 2006, Australia \\
\noindent Email: zhangou@maths.usyd.edu.au \\
\noindent Fax: + 61 2 9351 4534}

\end{document}